\documentclass [10pt,a4paper]{article}
\usepackage{amssymb}
\usepackage{enumerate}
\usepackage{amsmath}
\usepackage{amsfonts,mathrsfs}
\usepackage{amsthm}
\usepackage{epsfig}

\newcounter{teoremaganso}

\renewenvironment{abstract}{\small\quotation\noindent
 {\bfseries \abstractname .}}{\endquotation \par}

\newenvironment{trivlistparm}[2]{\trivlist
\item[{#1 #2}]}{\endtrivlist}

\newenvironment{prooftext}[1]{\trivlistparm{\bfseries}{#1}}{\Qed\endtrivlistparm}
\newenvironment{prova}{\trivlistparm{\bfseries}{Proof.}}{\Qed\endtrivlistparm}

\catcode`\@=11

\def\resetthefootnote{\renewcommand{\thefootnote}{\@arabic\c@footnote} }
\def\@principiremex#1{\trivlist
 \item[\hskip \labelsep{\bfseries #1\ \thetheo}]\ignorespaces}
\def\opar@principiremex#1[#2]{\trivlist
 \item[\hskip \labelsep{\bfseries #1\ \thetheo\ (#2)}]\ignorespaces}

\newcommand{\newTHEOremrom}[2]{\newenvironment{#1}{\refstepcounter{theo}\@ifnextchar[{\opar@principiremex{#2}}
{\@principiremex{#2}}}{\qedB\endtrivlist}} \catcode`\@=12
\DeclareMathSymbol{\square}{\mathord}{AMSa}{"03}
\newcommand{\qedB}{\nopagebreak\hspace*{\fill}$\square$\par}
\newcommand{\Qed}{\nopagebreak\hspace*{\fill}{\vrule width6pt height6pt depth0pt}\par}

\newtheorem {theo} {Theorem} [section]
\newtheorem {prop} [theo] {Proposition}

\newtheorem {lem} [theo] {Lemma}
\newtheorem {bigtheo} [teoremaganso] {Theorem}

\newTHEOremrom {defi} {Definition}
\newTHEOremrom {obs} {Remark}
\newTHEOremrom {ex} {Example}

\newcommand{\refc}[1]{\mbox{$(\ref{#1})$}}
\newcommand{\secc}[1]{Section~\ref{#1}}
\newcommand{\teoc}[1]{Theorem~\ref{#1}}
\newcommand{\propc}[1]{Proposition~\ref{#1}}

\newcommand{\lemc}[1]{Lemma~\ref{#1}}
\newcommand{\defic}[1]{Definition~\ref{#1}}
\newcommand{\obsc}[1]{Remark~\ref{#1}}

\newcommand{\figc}[1]{Figure~\ref{#1}}

\newcommand{\N}{\ensuremath{\mathbb{N}}}

\newcommand{\R}{\ensuremath{\mathbb{R}}}

\newcommand{\op}{\ensuremath{\mbox{\rm o}}}

\def\map#1#2#3{\mbox{${#1}\!:{#2}\longrightarrow{#3}$}}

\newcommand{\eps}{\varepsilon}

\newcommand{\sist}[2]{
  \left\{
   \begin{array}{l}
    \dot x=#1 \\[2pt] \dot y=#2
   \end{array}
  \right.
}

\begin{document}

\title{\textbf{Bifurcation of critical periods \\ from Pleshkan's isochrones}
 \footnotetext{{\it Key
words and phrases}: critical period; isochronous center; period
function; bifurcation; unfolding.}\footnotetext{M. Grau is partially
supported by the MCYT/FEDER grant number MTM2005-06098-C02-02. J.
Villadelprat is partially supported by the CONACIT through the grant
2005-SGR-00550 and by the DGES through the grant
MTM-2005-06098-C02-1.}}

\author{M. Grau and J. Villadelprat
\\*[.1truecm]
{\small \textsl{Departament de Matem\`{a}tica, }}
\\*[-.05truecm]
{\small \textsl{Universitat de Lleida, Lleida, Spain}}
\\*[.1truecm]
{\small \textsl{Departament d'Enginyeria Inform\`{a}tica i Matem\`{a}tiques,}}
\\*[-.05truecm]
{\small \textsl{Universitat Rovira i Virgili, Tarragona, Spain}}}

\date{}

\maketitle

\begin{abstract}
Pleshkan proved in 1969 that, up to a linear transformation and a
constant rescaling of time, there are four isochrones in the
family of cubic centers with homogeneous nonlinearities $\mathscr
C_3.$ In this paper we prove that if we perturb any of these
isochrones inside $\mathscr C_3,$ then at most two critical
periods bifurcate from its period annulus. Moreover we show that,
for each $k=0,1,2,$ there are perturbations giving rise to exactly
$k$ critical periods. As a byproduct, we obtain a partial result
for the analogous problem in the family of quadratic centers
$\mathscr C_2.$ Loud proved in 1964 that, up to a linear
transformation and a constant rescaling of time, there are four
isochrones in $\mathscr C_2.$ We prove that if we perturb three of
them inside $\mathscr C_2,$ then at most one critical period
bifurcates from its period annulus. In addition, for each $k=0,1,$
we show that there are perturbations giving rise to exactly $k$
critical periods. The quadratic isochronous center that we do not
consider displays some peculiarities that are discussed at the end
of the paper.
\end{abstract}

\section{Introduction and statement of the results \label{sect1}}

In this paper we study the period function of the centers of planar
vector fields with homogeneous nonlinearities, more concretely,
differential systems of the form
\begin{equation}\label{Hn}
\sist{-y+P_n(x,y),}{x+Q_n(x,y),}
\end{equation}
where $P_n(x,y)$ and $Q_n(x,y)$ are homogeneous real polynomials
of degree $n$. We denote by $\mathscr H_n$ the family of vector
fields of the above form and by $\mathscr C_n$ the subset of those
systems in $\mathscr H_n$ with a center at the origin. Finally,
$\mathcal I_n$ stands for the set of nonlinear isochrones in
$\mathscr H_n.$ Accordingly $\mathcal I_n\subset\mathscr
C_n\subset\mathscr H_n.$ We restrict ourselves to the cases $n=2$
and $n=3,$ that are the degrees for which the centers and the
isochrones have been classified.

Recall that the \emph{period annulus} of a center is the biggest
punctured neighbourhood foliated by periodic orbits and in what
follows we shall denote it by~$\mathcal P.$  Compactifying~$\R^2$
to the Poincar\'{e} disc, the boundary of~$\mathcal P$ has two
connected components: the center itself and a polycycle. We call
them respectively the \emph{inner} and \emph{outer boundary} of
the period annulus. The \emph{period function} of the center
assigns to each periodic orbit $\gamma$ in $\mathcal P$ its
period. If all the periodic orbits in $\mathcal P$ have the same
period, then the center is called \emph{isochronous.} Since the
period function is defined on the set of periodic orbits
in~$\mathcal P,$ usually the first step is to parameterize this
set, let us say $\{\gamma_s\}_{s\in (0,1)},$ and then one can
study the qualitative properties of the period function by means
of the map $s\longmapsto\mbox{period of $\gamma_s$},$ which is
analytic on~$(0,1).$ The \emph{critical periods} are the critical
points of this function and their number, character (maximum or
minimum) and distribution do not depend on the particular
parameterization of the set of periodic orbits used. We are
interested in the \emph{bifurcation} of critical periods. Roughly
speaking, the disappearance or emergence of critical periods as we
perturb a center. There are three different situations to study
(see~\cite{MMV2} for details):
\begin{enumerate}[$(a)$]
 \item Bifurcation of the period function from the inner boundary of $\mathcal
 P$ (i.e., the center itself).
 \item Bifurcation of the period function from $\mathcal P.$
 \item Bifurcation of the period function from the outer boundary of $\mathcal
 P$ (i.e., the polycycle).
\end{enumerate}

In this paper we are only concerned with bifurcation of the period
function from $\mathcal{P}$, i.e. case $(b)$, and when we perturb an
isochronous center. To be more precise, we study the bifurcation of
critical periods from $\mathcal P$ of the isochrones in $\mathcal
I_n$ when perturbed inside $\mathscr C_n$. Our initial motivation
was to study this problem for $n=3$ and the result we have obtained
is stated in \teoc{th1} below. As a byproduct of the tools we
develop to prove it, we get a partial result of the problem for
$n=2$, see \teoc{Loud}.

As the authors explain in \cite{Fran,FreGasGui1} many problems
dealing with the period function and critical periods are the
counterpart of problems about the Poincar\'{e} return map and
limit cycles. The problem studied here is the counterpart of the
question of how may limit cycles bifurcate from a center of an
integrable system when perturbed in a non-conservative direction.
The return map of the unperturbed system is the identity and one
asks for the number of fixed points that persist. In our case the
map that we perturb is constant (the period function of an
isochronous center) and we ask for the number of zeros of the
derivative that persist. Let us comment in this setting that the
bifurcation of the period function from $\mathcal P,$ i.e., case
$(b),$ when the unperturbed center is not isochronous is out of
reach for the moment. Its counterpart for limit cycles is the
so-called \emph{blue-sky} bifurcation, which is usually
undetectable.

The classification of $\mathscr C_3$ and $\mathcal I_3$ is due to
Malkin~\cite{Malkin} and Pleshkan~\cite{Pleshkan}, respectively. The
latter proved that, up to a linear change of coordinates and a
constant rescaling of time (see~\cite{ChaSab,MarMosRou}), there are
four isochronous centers in~$\mathscr H_3$:
\begin{equation}\label{llista}
\begin{array}{lll}
  &(S_1^{\ast})\quad \left\{
   \begin{array}{l}
    \dot x=-y -3x^2y+ y^3, \\ \dot{y}=x + x^3-3xy^2,
   \end{array}
  \right.\qquad
  &(S_2^{\ast})\quad   \left\{
   \begin{array}{l}
     \dot{x}=-y + x^2y,\\ \dot{y}=x+xy^2,
   \end{array}
  \right.\\[20pt]
   &(S_3^{\ast})\quad \left\{
   \begin{array}{l}
     \dot{x}=-y+3x^2y,\\ \dot{y}=x-2x^3+9xy^2,
   \end{array}
  \right.
  &(\bar S_3^{\ast})\quad   \left\{
   \begin{array}{l}
    \dot{x}=-y-3x^2y, \\ \dot{y}=x+2x^3-9xy^2.
   \end{array}
  \right.
\end{array}
\end{equation}
The first result that we shall prove is the following. Let us point
out that in its statement $\mathcal P$ refers to the period annulus
of the center at the origin.

\begin{bigtheo}\label{th1}
If we perturb $(S_1^{\ast}),$ $(S_2^{\ast}),$ $(S_3^{\ast})$ or
$(\bar S_3^{\ast})$ inside $\mathscr C_3,$ then at most two
critical periods bifurcate from $\mathcal P.$ Moreover, for each
$k=0,1,2,$ there are perturbations that give rise to exactly $k$
critical periods bifurcating from $\mathcal P.$
\end{bigtheo}

To be more precise, if~$X_0$ denotes the vector field associated to
one of the systems in \refc{llista}, then we consider the unfolding
given by $X_{\varepsilon}=X_0+Y_{\varepsilon}$ with
 \begin{equation}\label{pert1}
  Y_{\varepsilon}(x,y)\!:=\bigl(a(\varepsilon)x^2y+b(\varepsilon) y^3\bigr)\partial_x
     +\bigl(c(\varepsilon) x^3+d(\varepsilon)xy^2\bigr)\partial_y
 \end{equation}
where $a$, $b$, $c$ and $d$ are analytic functions vanishing at
$\varepsilon=0.$ Since $X_{\eps}$ is reversible, it is clear that
it belongs to $\mathscr C_3$ for all $\eps.$ The fact that there
is no other perturbation with this property follows from Malkin's
result (see also \cite{RouToni}). Of course the problem under
consideration only makes sense in case that the perturbation is
not isochronous itself. This, let us say, \emph{isochronous
direction of perturbation} that we must avoid is $d+3c=a+3c=b-c=0$
for $(S_1^{\ast});$ $b=c=a-d=0$ for $(S_2^{\ast});$ and
$b=2a+3c=2d+9c=0$ for $(S_3^{\ast})$ and $(\bar S_3^{\ast}).$ As
we will see, the existence of such a direction prevents the
parameters in the perturbation from acting independently and it
explains the reason why three critical periods cannot appear even
though we consider a 4-parameter unfolding.

There are some previous results related to \teoc{th1} that should be
referred. Rousseau and Toni prove in~\cite{RouToni} that if we
perturb any isochronous center in \refc{llista} inside $\mathscr
C_3,$ then at most two critical periods bifurcate from the inner
boundary of $\mathcal P.$ They also study the perturbation of the
linear (isochronous) center $X_0=-y\partial_x+x\partial_y.$ In this
case there are three different unfoldings of $X_0$ inside $\mathscr
C_3$ and the authors prove that they give rise, respectively, to at
most $0$, $1$ and $3$ critical periods bifurcating from the inner
boundary of $\mathcal P$. In \cite{CimGasSil} the authors also study
unfoldings of the linear center. They provide some lower bounds for
the maximum number of critical periods that can bifurcate from
$\mathcal P.$

The classification of $\mathscr C_2$ and $\mathcal I_2$ is due to
Dulac~\cite{Dulac} and Loud~\cite{Loud}, respectively. The latter
proved that, up to a linear change of coordinates and a constant
rescaling of time, there are four isochronous centers in~$\mathscr
H_2$:
\begin{equation}\label{llista2}
\begin{array}{lll}
  &(S_1)\quad \sist{-y+xy,}{x-\frac{1}{2}x^2+\frac{1}{2}y^2,}\qquad
  &(S_2)\quad \sist{-y+xy,}{x+y^2,}\\[20pt]
  &(S_3)\quad \sist{-y+xy,}{x+\frac{1}{4}y^2,}
  &(S_4)\quad \sist{-y+xy,}{x-\frac{1}{2}x^2+2y^2.}
\end{array}
\end{equation}
Concerning these systems and the period annulus $\mathcal P$ of
the isochronous center at the origin, we prove the following
result:
\begin{bigtheo}\label{Loud}
If we perturb $(S_2),$ $(S_3)$ or $(S_4)$ inside $\mathscr C_2,$
then at most one critical period bifurcates from $\mathcal P.$
Moreover, for each $k=0,1,$ there are perturbations that give rise
to exactly $k$ critical periods bifurcating from $\mathcal P.$
\end{bigtheo}

The perturbation of the isochrones $(S_2),$ $(S_3)$ and $(S_4)$
inside $\mathscr C_2$ can be studied with the same unfolding.
Indeed, if $X_0$ is the associated vector field to one of these
three systems, then it suffices to consider
$X_{\eps}=X_0+Y_{\eps}$ with
 \begin{equation}\label{pert2}
  Y_{\varepsilon}(x,y)\!:=b(\varepsilon)xy\partial_x
     +\bigl(d(\varepsilon) x^2+f(\varepsilon)y^2\bigr)\partial_y
 \end{equation}
where $b$, $d$ and $f$ are analytic functions vanishing at
$\varepsilon=0.$ Clearly $X_{\eps}$ has a center at the origin
because the vector field is reversible. The fact that there is not
any other perturbation of~$X_0$ inside~$\mathscr C_2$ follows from
Dulac's result. On the contrary, to study the perturbation
of~$(S_1)$ inside~$\mathscr C_2$ it is necessary to consider two
unfoldings. The first one corresponds to the perturbation inside
the reversible centers and the second one to the perturbation
inside the so-called generalized Lotka-Volterra centers
(see~\cite{ChiJac} for details). However this is not the reason
why \teoc{Loud} does not contemplate the perturbation of~$(S_1).$
We explain what happens in this case at the end of the paper.

There are previous results related with \teoc{Loud} that must be
referred. Chicone and Jacobs~\cite{ChiJac} prove that if we perturb
any isochronous center in \refc{llista2} inside $\mathscr C_2,$ then
at most one critical period bifurcates from the inner boundary of
$\mathcal P$ (i.e., the center). They consider as well the
perturbation of the linear isochronous center, for which the upper
bound is two. The bifurcation of critical periods from the outer
boundary of $\mathcal P$ (i.e, the polycycle) is partially studied
in \cite{MMV1,MMV2,MV}. (In \cite{MMV2} the authors present a
conjectural bifurcation diagram of the period function of any
quadratic center in the Loud form. We reproduce this diagram at the
end of the paper, see \figc{dib}, where we shall comment how
\teoc{Loud} and this conjecture fit together.) Gasull and Yu
\cite{GasYu} prove \teoc{Loud} assuming that the perturbation
\refc{pert2} is \emph{linear} in~$\varepsilon$ (or, as they say,
perturbations up to first order in $\eps$). Let us note moreover
that system $(S_4)$ has another isochronous center apart from the
one at the origin and that the authors in \cite{GasYu} study the
simultaneous bifurcation of critical periods from both period
annuli. In this paper we restrict our study to the bifurcation of
critical periods from the period annulus of the center at the
origin. 
Finally, the assertion in \teoc{Loud} corresponding to the
isochronous center $(S_2)$ is proved in~\cite{GasZhao}.

In view of Theorems~\ref{th1} and \ref{Loud}, and all the previous
results we refer above, one may pose the following question for the
family~$\mathscr H_n$ of vector fields with homogeneous
nonlinearities of degree~$n$: Is it true that at most $n-1$ critical
periods can bifurcate from a nonlinear isochronous center when
perturbed inside~$\mathscr C_n?$

The proof of \teoc{th1} is given in \secc{sect3} and it is based on
a recent result which appears in~\cite{GMV}, see the forthcoming
Proposition \ref{prop2}. This recent result is a criterion for a
collection of Abelian integrals to have a Chebyshev property and it
provides, in many cases, a shortcut with respect to the classical
methods for the same issue. The proof for system $(S_1^{\ast})$
follows the classical approach and, thus, it is much longer than the
proof for the other three cases. The proof of \teoc{Loud} is given
in \secc{sect4}. The following section contains the definitions and
previous results needed to prove our results.

\section{Definitions and notation \label{sect2}}

Let $\{X_{\varepsilon}\}$ be an analytic family of planar vector
fields with a center at the origin for all $\varepsilon\approx 0$
and assume that the center is isochronous for $\varepsilon=0.$ Let
$\mathcal P$ be the \emph{period annulus} of the center of $X_0$
(i.e., the biggest punctured neighbourhood of the origin foliated by
periodic orbits) and let \map{\xi}{\mathcal I}{\R^2} be an analytic
transversal section to $X_0$ in $\mathcal P$, with $\mathcal{I}$ an
open real interval. Let $T(s;\varepsilon)$ be the period of the
periodic orbit of $X_{\varepsilon}$ passing through the point
$\xi(s)\in\mathcal P$ whenever is defined. (Here we use that the
curve given by the graphic of $\xi$ is still a transversal section
to~$X_{\varepsilon}$ for $\varepsilon\approx 0.)$ We can choose, for
instance, the transversal section to be a solution of the vector
field orthogonal to $X_0$ and in this case $\mathcal I$ is its
maximal interval of definition.

\begin{defi}\label{defi2}
We say that $k$ \emph{critical periods bifurcate from $\mathcal
P$} in the family $\{X_{\varepsilon}\}$ as $\varepsilon
\rightarrow 0$ if there exist $k$ functions $s_i(\varepsilon),$
continuous in a neighbourhood of $\varepsilon=0,$ and such that
$T'\bigl(s_i(\varepsilon);\varepsilon\bigr)\equiv 0$ and
$s_{i}(0)=s_i^{\star}\in\mathcal I$ for each $i=1,2,\ldots,k.$
\end{defi}

Since $T(s;\varepsilon)$ is analytic for $\varepsilon\approx 0,$ we
can consider its Taylor development at $\varepsilon=0,$ say
 \[
  T(s;\varepsilon)=\sum_{i=0}^{\infty}T_i(s)\varepsilon^i.
 \]
Note that in fact $T_0$ is constant because by assumption the
center is isochronous for $\varepsilon=0$. Then, if the center
of~$X_{\varepsilon}$ is non-isochronous for $\varepsilon\neq 0,$
there exists some $\ell\geqslant 1$ such that
 \[
  T'(s;\varepsilon)=T_{\ell}'(s)\varepsilon^{\ell}+\op(\varepsilon^{\ell}),
 \]
where $T_{\ell}'$ is not identically zero and the remainder is
uniform in $s$ on each compact subinterval of $\mathcal I.$ In this
case, from
\defic{defi2} it follows that $T'_{\ell}(s_{i}^{\star})$ for
$i=1,2,\ldots,k.$ Consequently, by applying the Weierstrass
Preparation Theorem, the number of zeros of $T_{\ell}'(s)$ for $s\in
\mathcal I,$ counted with multiplicities, provides an upper bound
for $k.$ A lower bound for $k$ is given by the number of simple
zeros of $T_{\ell}'(s)$ in $\mathcal I$ by using the Implicit
Function Theorem.

\begin{obs}\label{indep}
The number of critical periods that bifurcate from $\mathcal P$ as
$\varepsilon\rightarrow 0$ does not depend on the particular
parameterization of the set of periodic orbits in $\mathcal P$
used. If two of these parameterizations yield to
$T(s;\varepsilon)$ and $\widehat T(s;\varepsilon),$ then there
exists a diffeomorphism \map{\zeta}{\widehat {\mathcal
I}}{\mathcal I} such that $T(\zeta(s);\varepsilon)=\widehat
T(s;\varepsilon).$ Following the obvious notation, this implies
that $\widehat
T_{\ell}'(s)=\zeta'(s)T_{\ell}'\bigl(\zeta(s)\bigr).$
\end{obs}

\begin{defi}\label{transversal}
On account of the above observation, from now on we shall work
with a transversal section constructed in the following way. We
take a vector field $U_0$ analytic on $\mathcal P\cup\{(0,0)\},$
transversal to~$X_0$ at~$\mathcal P$ and such that
$[X_0,U_0]\equiv 0.$ In the literature such a vector field is
called a \emph{commutator} of~$X_0$ and it always exists in a
neighbourhood of an isochronous center, see \cite{Sabatini}. Here
we require more; it must be analytic on $\mathcal P\cup\{(0,0)\},$
but in all the examples that we shall study such a commutator
exists. Then we define $T(s,\varepsilon)$ exactly as before taking
a solution of $U_0$ as the transversal section~$\xi.$ More
concretely, we choose any $q\in\mathcal P$ and we consider the
solution $\psi(s;q)$ of~$U_0$ with $\psi(0;q)=q.$ Then we define
$T(s;\varepsilon)$ as the period of the periodic orbit of
$X_{\varepsilon}$ passing through $\xi(s)\!:=\psi(s;q).$
\end{defi}

Given an analytic family of functions $\{g_{\varepsilon}\}$ we write
its Taylor series at $\varepsilon=0$ as
$g(x,y)=\sum_{i=0}^{\infty}g_i(x,y)\varepsilon^i.$ We shall also use
the notation $j^{\,k}(g)=\sum_{i=0}^kg_i\varepsilon^i$ for its
$k$-jet with respect to $\varepsilon$. (This notation extends for
vector fields in the obvious way.) Moreover, in what follows, we
shall say that two vector fields $X$ and $Y$ are \emph{conjugated}
if there exists a change of coordinates $\Phi$ transforming $X$ to
$Y$, i.e., $\Phi^{\ast}X=Y$, where
\[
\bigl(\Phi^{\ast}X\bigr)(p)=\bigl(D\Phi\bigr)^{-1}_p\bigl(X\circ\Phi(p)\bigr).
\]
Then we say that $\Phi$ \emph{conjugates} $X$ and $Y.$ A classical
result shows that $X$ has an isochronous center at $p_0$ if and
only if there exists an analytic diffeomorphism $\Phi$ on a
neighbourhood of $p_0$ that conjugates $X$ with its linear part at
$p_0.$ In this case $\Phi$ is called a $\emph{linearization}.$

\section{Proof of Theorem \ref{th1} \label{sect3}}

Let $X_0$ denote the vector field associated to one of the systems
listed in (\ref{llista}) and $\mathcal{P}$ be the period annulus
of the isochronous center at the origin. Let $U_0$ be a commutator
of $X_0$ which is analytic on $\mathcal P\cup\{(0,0)\}$. We
describe each of the considered vector fields $U_0$ in Lemma
\ref{lem1} and we note that they are, indeed, polynomial vector
fields. Since $X_0$ and $U_0$ are transversal vector fields on
$\mathcal P$, there exist two analytic functions $\lambda$ and
$\mu$ such that
 \[
  Y_{\eps}=X_{\varepsilon}-X_0\, =\, \lambda\, X_0\, +\, \mu\, U_0,
 \]
where $Y_{\varepsilon}$ is the vector field corresponding to the
perturbative part as defined in (\ref{pert1}). In fact it is clear
that
\[
 \lambda\, =\, \frac{Y_{\eps}\wedge U_0}{X_0\wedge U_0}\quad \mbox{ and }
 \quad \mu\, =\, \frac{X_0\wedge Y_{\eps}}{X_0\wedge U_0}.
\]
Our first result, \teoc{teo1}, gives an expression for the
$(k+1)$th term of the Taylor series of $T'(s;\varepsilon)$ at
$\varepsilon=0.$ It depends on the $(k+1)$th term of the
development of $\lambda=\lambda(x,y;\varepsilon)$ at
$\varepsilon=0$ and to prove it we shall use the following result
(see~\cite{GasYu}). This theorem is a generalization of the result
stated in the seminal paper \cite{FreGasGui2}, see also the
references therein.

\begin{theo}[Gasull-Yu]\label{Gasull-Yu}
Let $X$ be an analytic vector field with a center at $p$ and let
$\mathcal P$ be its period annulus. Let $U$ be an analytic vector
field on $\mathcal P\cup\{p\}$ which is transversal to $X$ on
$\mathcal P$. Then $[X,U]=\alpha X+\beta U$ for some analytic
functions $\alpha$ and $\beta.$ Fix any $q\in\mathcal P$ and let
\map{\xi}{\mathcal I}{\mathcal P} be the solution of $U$ with
$\xi(0)=q.$ For each $s\in\mathcal I,$ let $\varphi(t;s)$ be the
solution of $X$ with $\varphi(0;s)=\xi(s)$ and let $T(s)$ be its
period. Then
 \[
  T'(s)=\int_0^{T(s)}\alpha\bigl(\varphi(t;s)\bigr)
        \exp\left(-\int_0^t\beta\bigl(\varphi(\tau;s)\bigr)d\tau\right)dt\,
  \mbox{ for all $s\in\mathcal I.$}
 \]
\end{theo}

Let us point out that in the above statement $U$ is just a
transversal vector field to $X$. We can now prove the following:

\begin{theo}\label{teo1}
Consider the parameterization of the period function $T(s;\eps)$
as introduced in \defic{transversal}. Let us assume that, for some
$k\in\N,$ there exists an analytic family of
diffeomorphisms~$\{\Phi_{\varepsilon}\},$ in a neighbourhood of
$(0,0)$, such that $\Phi_{\varepsilon}$ linearizes
$j^k\bigl(X_{\varepsilon}\bigr)$ for each $\varepsilon\approx 0.$
Then $T_0'\equiv T_1'\equiv\ldots\equiv T_k'\equiv 0$ and
 \[
  T_{k+1}'(s)=-\int_0^{T_0}\left.U_0(\lambda_{k+1})\right|_{(x,y)=\varphi(t;s)}dt\,
  \mbox{ for all $s\in\mathcal I,$}
 \]
where $\varphi(t;s)$ is the solution of $X_0$ with
$\varphi(0;s)=\xi(s).$
\end{theo}

\begin{prova}
We focus in the formula for $T_{k+1}'$ because as a byproduct it
will follow that $T_0'\equiv T_1'\equiv\ldots\equiv T_k'\equiv 0$.
Note first that it suffices to show the validity of the equality for
periodic orbits near the origin since the functions on the left and
right of the equality depend analytically on~$s.$ To this end we use
the family of diffeomorphisms $\{\Phi_{\varepsilon}\}$ that are
defined in neighbourhood of the origin. (Certainly this family
depends on $k$ as well, but for the sake of shortness we omit this
additional subscript.)

Define
$\Psi_{\varepsilon}\!:=(\Phi_0)^{-1}\!\circ\Phi_{\varepsilon}.$ Then
$\{\Psi_{\varepsilon}\}$ is an analytic family of diffeomorphisms
with $\Psi_{\varepsilon}=\mbox{Id}+\op(\varepsilon^0)$ and such
that, for each $\varepsilon\approx 0,$ $\Psi_{\varepsilon}$
conjugates $j^{k}(X_{\varepsilon})$ and $X_0,$ i.e.,
$(\Psi_{\varepsilon})^{\ast}\bigl(j^k(X_{\varepsilon})\bigr)=X_0.$
Let us define in addition
$U_{\varepsilon}\!:=\left(\Psi_{\varepsilon}\right)_{\ast}(U_0)$ and
note that then $U_{\varepsilon}=U_0+\op(\varepsilon^0).$ We claim
that
 \begin{equation}\label{eq2}
  [X_{\varepsilon},U_{\varepsilon}]=
  \varepsilon^{k+1}\bigl(-U_0(\lambda_{k+1})+\op(\varepsilon^0)\bigr)X_{\varepsilon}+
  \varepsilon^{k+1}\bigl(-U_0(\mu_{k+1})+\op(\varepsilon^0)\bigr)U_{\varepsilon},
 \end{equation}
where here (and in what follows) we use the notation
$\lambda=\sum_{i=0}^{\infty}\lambda_i\varepsilon^i$ and
$\mu=\sum_{i=0}^{\infty}\mu_i\varepsilon^i.$ In order to prove the
claim let us first recall that $X_{\varepsilon}=X_0+\lambda
X_0+\mu U_0.$ Thus, we can write
 \begin{equation}\label{eq2.5}
  X_{\varepsilon}=j^k(X_{\varepsilon})+\varepsilon^{k+1}\Bigl(\bigl(\lambda_{k+1}+\op(\varepsilon^0)\bigr)X_0+
  \bigl(\mu_{k+1}+\op(\varepsilon^0)\bigr)U_0\Bigr).
 \end{equation}
On the other hand,
 \[
 (\Psi_{\varepsilon})^{\ast}\!\left[j^k(X_{\varepsilon}),U_{\varepsilon}\right]=
  \left[(\Psi_{\varepsilon})^{\ast}\bigl(j^k(X_{\varepsilon})\bigr),(\Psi_{\varepsilon})^{\ast}(U_{\varepsilon})\right]
  \\
  =[X_0,(\Psi_{\varepsilon})^{\ast}\bigl(\left(\Psi_{\varepsilon}\right)_{\ast}(U_0)\bigr)]=
  [X_0,U_0]=0.
 \]
In the first equality above we use that if $f$ is a local
diffeomorphism, then $f^{\ast}[X,Y]=[f^{\ast}X,f^{\ast}Y],$ see
for instance \cite{KMS}. In the second one we use that
$\Psi_{\varepsilon}$ conjugates $j^{k}(X_{\varepsilon})$ and $X_0$
and that, by definition, $U_{\varepsilon} = \left(
\Psi_{\varepsilon} \right)_{\ast} (U_0)$. The last one follows
from the fact that $U_0$ is a commutator of $X_0$. We can thus
conclude that $[j^k(X_{\varepsilon}),U_{\varepsilon}]=0.$
Accordingly, from \refc{eq2.5} we obtain
 \[
  [X_{\varepsilon},U_{\varepsilon}]=\underbrace{[j^k(X_{\varepsilon}),U_{\varepsilon}]}_{0}+
  \varepsilon^{k+1}\underbrace{\left[\bigl(\lambda_{k+1}+\op(\varepsilon^0)\bigr)X_0,U_{\varepsilon}\right]}_{\Delta_1}+
  \varepsilon^{k+1}\underbrace{\left[\bigl(\mu_{k+1}+\op(\varepsilon^0)\bigr)U_0,U_{\varepsilon}\right]}_{\Delta_2}.
 \]
Then, since $[fX,Y]=f[X,Y]-Y(f)X$ (see \cite{KMS} again), it follows
that
 \begin{align*}
  \Delta_1&=\bigl(\lambda_{k+1}+\op(\varepsilon^0)\bigr)[X_0,U_{\varepsilon}]
  -U_{\varepsilon}\bigl(\lambda_{k+1}+\op(\varepsilon^0)\bigr)X_0\\[3pt]
  &=\bigl(\lambda_{k+1}+\op(\varepsilon^0)\bigr)[X_0,U_0+\op(\varepsilon^0)]
  -U_0(\lambda_{k+1}\bigr)X_{\varepsilon}+\op(\varepsilon^0)
  =-U_0(\lambda_{k+1})X_{\varepsilon}+\op(\varepsilon^0)
  \intertext{and}
  \Delta_2&=\bigl(\mu_{k+1}+\op(\varepsilon^0)\bigr)[U_0,U_{\varepsilon}]
  -U_{\varepsilon}\bigl(\mu_{k+1}+\op(\varepsilon^0)\bigr)U_0
  =-U_0(\mu_{k+1})U_{\varepsilon}+\op(\varepsilon^0),
 \end{align*}
where we took
$X_{\varepsilon}-X_0=U_{\varepsilon}-U_0=\op(\varepsilon^0)$ and
$[X_0,U_0]=0$ into account. Hence
 \[
  [X_{\varepsilon},U_{\varepsilon}]=
  \varepsilon^{k+1}\bigl(-U_0(\lambda_{k+1})\bigr)X_{\varepsilon}+
  \varepsilon^{k+1}\bigl(-U_0(\mu_{k+1})\bigr)U_{\varepsilon}+\op(\varepsilon^{k+1}),
 \]
and \refc{eq2} follows after decomposing the remainder term above in
the $X_{\varepsilon}$ and $U_{\varepsilon}$ components using that
both vector fields are transversal for $\varepsilon\approx 0.$ This
shows that the claim is true.

Next we take the point $q\in\mathcal P$ in \defic{transversal} and
consider the solution $\psi_{\varepsilon}(s;q)$ of $X_{\varepsilon}$
with $\psi_{\varepsilon}(0;q)=q.$ If $\widehat T(s;\varepsilon)$
denotes the period of the periodic orbit of $X_{\varepsilon}$
passing through $\xi_{\varepsilon}(s)\!:=\psi_{\varepsilon}(s;q),$
then from \refc{eq2} and applying \teoc{Gasull-Yu} it follows that
 \[
  \widehat T'(s;\varepsilon)=-\varepsilon^{k+1}\int_0^{\widehat T(s;\varepsilon)}
    \left.\!\bigl(U_0(\lambda_{k+1})+\op(\varepsilon^0)\bigr)\right|_{(x,y)=\varphi_{\varepsilon}(t;s)}
     e^{-\varepsilon^{k+1}\int_0^t\left.\left(
     U_{\varepsilon}(\mu_{k+1})+\op(\varepsilon^0)
     \right)\right|_{(x,y)=\varphi_{\varepsilon}(\tau;s)}d\tau}dt,
 \]
where $\varphi_{\varepsilon}(t;s)$ is the solution of
$X_{\varepsilon}$ with
$\varphi_{\varepsilon}(0;s)=\xi_{\varepsilon}(s).$ Therefore, on
account of
$X_{\varepsilon}-X_0=U_{\varepsilon}-U_0=\op(\varepsilon^0),$ we
obtain
 \[
  \widehat T'(s;\varepsilon)=-\varepsilon^{k+1}\int_0^{T_0}
    \left.\!U_0(\lambda_{k+1})\right|_{(x,y)=\varphi_0(t;s)}dt+\op(\varepsilon^{k+1}).
 \]
Note at this point that, by construction, $T(s;\varepsilon)=\widehat
T\bigl(f(s,\varepsilon);\varepsilon\bigr)$ for some function $f$
with $f(s,0)=s.$ Taking this into account, one can easily show from
the above equality that
 \[
   T'(s;\varepsilon)=-\varepsilon^{k+1}\int_0^{T_0}
    \left.\!U_0(\lambda_{k+1})\right|_{(x,y)=\varphi_0(t;s)}dt+\op(\varepsilon^{k+1}).
 \]
Accordingly, $T_0'\equiv T_1'\equiv\ldots\equiv T_k'\equiv 0$ and
 \[
  T_{k+1}'(s)=-\int_0^{T_0}\left.U_0(\lambda_{k+1})\right|_{(x,y)=\varphi(t;s)}dt.
 \]
This completes the proof of the result.
\end{prova}

\begin{obs} \label{cargol}
As we explain at the end of \secc{sect2}, it is well known that a
vector field has an isochronous center if and only if it is
linearizable. The authors do not know weather this
characterization is true for families. This is the reason why the
assumption in \teoc{teo1} requires the existence of an
\emph{analytic family} of linearizations instead of assuming that
the center of $j^k\bigl(X_{\varepsilon}\bigr)$ is isochronous for
$\varepsilon\approx 0$. The latter assumption provides a
linearization $\Phi_{\varepsilon}$ for each $\varepsilon\approx
0,$ but in general we know nothing about the regularity of the
family $\{\Phi_{\eps}\}$ with respect to $\eps.$
\end{obs}

The following lemma (see \cite{ChaSab, MarMosRou}) provides the
necessary information in order to apply the above result to our
problem with the Pleshkan's isochrones. It gives the commutator
and first integral of each isochronous center in \refc{llista}.

\begin{lem}\label{lem1}

\begin{itemize}
\item[$(S_1^{\ast})$] The vector field $X_0=(-y-3x^2y + y^3)
\partial_x +(x+x^3-3xy^2)\partial_y$ has first integral
\[
H(x,y)=\frac{x^2+y^2}{\sqrt{(x^2+(y-1)^2)(x^2+(y+1)^2)}}
\]
and commutator $U_0 =(x+x^3-3xy^2)\partial_x
+(y+3x^2y-y^3)\partial_y$.

\item[$(S_2^{\ast})$] The vector field $X_0=\bigl(-y + x^2y\bigr)
\partial_x +\bigl(x+xy^2\bigr)\partial_y$ has first integral
\[
H(x,y)=\frac{x^2+y^2}{1-x^2}
\]
and commutator  $U_0 = x(1-x^2)\partial_x + y(1-x^2)\partial_y$.

\item[$(S_3^{\ast})$] The vector field $X_0= (-y + 3x^2y)
\partial_x + (x-2x^3+9xy^2)\partial_y$ has first integral
\[
H(x,y)=\frac{(x-2x^3)^2+y^2}{(1-3x^2)^3}
\]
and commutator $U_0 =x(1-3x^2)(1-2x^2)\partial_x +
y(1-3x^2)(1-6x^2)\partial_y$.

\item[$(\bar{S_3^{\ast}})$] The vector field $X_0= (-y - 3x^2y)
\partial_x + (x+2x^3-9xy^2)\partial_y$ has first integral
\[
H(x,y)=\frac{(x+2x^3)^2+y^2}{(1+3x^2)^3}
\]
and commutator $U_0 =x(1+3x^2)(1+2x^2)\partial_x +
y(1+3x^2)(1+6x^2)\partial_y$.

\end{itemize}
\end{lem}

Since the number of critical periods that bifurcate from $\mathcal
P$ does not depend on the particular parameterization of the set
of periodic orbits in $\mathcal P$ used (recall \obsc{indep}),
\teoc{th1} will follow once we prove the next result.

\begin{theo}\label{teo2}
Let $X_0$ be one of the four Pleshkan's isochrones in \lemc{lem1}
and let $U_0$ be its corresponding commutator. Consider the
unfolding of centers $X_{\varepsilon}\!:=X_0+ Y_{\varepsilon}$ where
 \[
  Y_{\varepsilon}=\bigl(a(\varepsilon)x^2y+b(\varepsilon)y^3\bigr)\partial_x
     +\bigl(c(\varepsilon)x^3+d(\varepsilon)xy^2\bigr)\partial_y
 \]
and where $a$, $b$, $c$ and $d$ are analytic functions in a
neighbourhood of $\varepsilon=0$ with $a(0)=b(0)=c(0)=d(0)=0.$ Let
\map{\xi}{\mathcal I}{\R^2} be a transversal section to $\mathcal
P$ given by a solution of $U_0$ and let $T(s;\varepsilon)$ be the
period of the periodic orbit of $X_{\varepsilon}$ passing through
$\xi(s)\in\mathcal P.$ If $T_0'\equiv T_1'\equiv\ldots T_k'\equiv
0$ and $T_{k+1}'\not\equiv 0,$ then $T_{k+1}'(s)$ has at most two
zeros for $s\in\mathcal I$ and there are perturbations
$Y_{\varepsilon}$ giving rise to $0,1$ and $2$ zeros.
\end{theo}

In order to prove this result some definitions and lemmas are
needed.

\begin{defi}\label{defi1}
Let $f_0,f_1,\ldots,f_{n-1}$ be analytic functions on an open
interval $L$ of $\R.$
\begin{enumerate}[$(a)$]
\item $(f_0,f_1,\ldots,f_{n-1})$ is a \emph{Chebyshev system} $($in short, T-system$)$ on
      $L$ if any nontrivial linear combination
      \[
       \alpha_0f_0(x)+\alpha_1f_1(x)+\ldots+\alpha_{n-1}f_{n-1}(x)
      \]
      has at most $n-1$ isolated zeros on $L.$
\item $(f_0,f_1,\ldots,f_{n-1})$ is a \emph{complete Chebyshev
      system} $($in short, CT-system$)$ on $L$ if
      $(f_0,f_1,\ldots,f_{k-1})$ is a T-system for all $k=1,2,\ldots,n.$
\item $(f_0,f_1,\ldots,f_{n-1})$ is an \emph{extended complete
      Chebyshev system} $($in short, ECT-system$)$ on $L$ if, for all $k=1,2,\ldots,n,$ any
      nontrivial linear combination
      \[
       \alpha_0f_0(x)+\alpha_1f_1(x)+\ldots+\alpha_{k-1}f_{k-1}(x)
      \]
      has at most $k-1$ isolated zeros on $L$ counted with multiplicities.
\end{enumerate}
(Let us mention that, in these abbreviations, ``T'' stands for
Tchebycheff, which in some sources is the transcription of the
Russian name Chebyshev.)
\end{defi}

\begin{obs}\label{hison}
If $(f_0,f_1,\ldots,f_{n-1})$ is an ECT-system on $L,$ then for each
$k=0,1,\ldots,n-1$ there exists a linear combination
      \[
       \alpha_0f_0(x)+\alpha_1f_1(x)+\ldots+\alpha_{n-1}f_{n-1}(x)
      \]
with \emph{exactly} $k$ simple zeros on $L$ (see for instance
\cite{Karlin,Mar}).
\end{obs}

\begin{defi}
Let $f_0,f_1,\ldots,f_{k-1}$ be analytic functions on an open interval $L$ of $\R.$ The
\emph{Wronskian} of $(f_0,f_1,\ldots,f_{k-1})$ at $x\in L$ is
 \[
 W\bigl[f_0,f_1,\cdots,f_{k-1}\bigr](x)=\det\left(f_j^{(i)}(x)\right)_{0\leqslant i,j\leqslant k-1}
   =\left|\begin{array}{ccc}
    f_0(x) & \cdots & f_{k-1} (x) \\
    f'_0(x) & \cdots & f'_{k-1} (x) \\
    & \vdots & \\
    f^{(k-1)}_{0}(x) & \cdots & f^{(k-1)}_{k-1} (x) \\
    \end{array} \right|.
  \]
\end{defi}

The following result is well known (see \cite{Karlin,Mar} for instance).

\begin{lem}\label{lem2}
$(f_0,f_1,\ldots,f_{n-1})$ is an ECT-system on $L$ if, and only if, for each
$k=1,2,\ldots,n$ it holds
 \[
  W\bigl[f_0,f_1,\cdots,f_{k-1}\bigr](x)\neq 0\mbox{ for all $x\in L$}.
 \]
\end{lem}

Recall that \teoc{teo2} deals with four different isochrones. For
each one, the proof yields to a linear combination of some Abelian
integrals and it is necessary to show that they form an ECT-system.
This is in general extremely complicated to verify, but in three of
the cases it is not that difficult because a criterion proved in
\cite{GMV} applies successfully. \propc{prop2} is a simplified
version of this criterion that it suffices for our purposes. In its
statement we suppose that the projection of $\mathcal P$ on the
$x$-axis is given by $(-x_r,x_r)$ and that the energy level of~$H$
at the periodic orbits in~$\mathcal P$ ranges from $h=0$ to $h=h_0.$

\begin{prop}\label{prop2}
Let us consider the Abelian integrals
\[
I_i(h)=\int_{\gamma_h}f_i(x){y^{2m-1}}dx,\ \mbox{
$i=0,1,\ldots,n-1,$}
\]
where $f_i(x)$ are analytic functions in $(-x_r,x_r)$, $m \in
\mathbb{Z}$ and where, for each $h\in (0,h_0),$ $\gamma_h$ is the
oval surrounding the origin inside the level curve
$\{A(x)+B(x)y^{2}=h\}.$ Assume that $A$ and $B$ are even functions
and let $\ell_i$ be the even part of $f_i.$ Then
$(I_0,I_1,\ldots,I_{n-1})$ is an ECT-system on $(0,h_0)$ if
$m\geqslant n-1$ and $\bigl(\ell_0,\ell_1,\ldots,\ell_{n-1}\bigr)$
is a CT-system on $(0,x_r).$
\end{prop}

The following lemma, proved in \cite{GMV}, establishes a formula to
write the integrand of an Abelian integral so as to be suitable to
apply the \propc{prop2}.

\begin{lem}\label{puja} Let $\gamma_h$ be the
oval surrounding the origin inside the level curve
$\{A(x)+B(x)y^2=h\}$ and we consider a function $F$ such that
$F/A'$ is analytic at $x=0.$ Then, for any $k\in\N,$
 \[
  \int_{\gamma_h}F(x)y^{k-2}dx=\int_{\gamma_h} G(x)y^kdx
 \]
where $G(x)=\frac{2}{k}\left(\frac{BF}{A'}\right)'(x)-
\left(\frac{B'F}{A'}\right)(x).$
\end{lem}

The criterion in \propc{prop2} does not apply for the isochronous
center $(S_1^{\ast}).$ This fact makes the proof of \teoc{teo2} for
this case much longer than the others. In particular we shall need
some properties of the complete elliptic integrals (see \cite{Abram}
for instance).

\begin{defi}\label{complete}
The complete elliptic integrals of first and second kind are
 \[
  \mathcal{E}(u)\!:=\int_{0}^{\pi/2}\sqrt{1-u^2 \sin^2 t}\,dt\,\mbox{ and }\,
  \mathcal{K}(u)\!:=\int_{0}^{\pi/2}\frac{dt}{\sqrt{1-u^2 \sin^2 t}}
 \]
respectively, which are analytic real functions for $u\in (-1,1).$
\end{defi}

\begin{lem}\label{lem5} The functions $\mathcal K$ and $\mathcal E$ verify the linear differential
equation
 \[
  \left(\!\begin{array}{c} \mathcal K' \\ \mathcal E'\end{array}\!\right)=\frac{-1}{u}
  \left(\begin{array}{cc} 1 & \frac{1}{u^2-1} \\ 1 & -1 \end{array}\right)
  \left(\!\begin{array}{c} \mathcal K \\ \mathcal E \end{array}\!\right).
 \]
Moreover, their Taylor series at $u=0$ are
\[
\mathcal{K}(u)=\frac{\pi}{2} \, \sum_{i\geqslant0} \left(
\frac{(2i-1)!!}{(2i)!!} \right)^2 u^{2i}\,\mbox{ and
}\,\mathcal{E}(u) \, = \, -\, \frac{\pi}{2} \, \sum_{i\geqslant0}
\left( \frac{(2i-1)!!}{(2i)!!} \right)^2 \frac{u^{2i}}{2i-1},
\]
where $u\in (-1,1).$
\end{lem}

Recall that the double factorial of an integer $n$, with $n\geqslant
-1$, is defined as
\[
n!!=\left\{ \begin{array}{ll} n(n-2) \ldots 5 \cdot 3 \cdot 1 & \ \
\mbox{if $n>0$ and $n$ odd}, \\ n(n-2) \ldots
6 \cdot 4 \cdot 2 & \ \ \mbox{if $n>0$ and $n$ even}, \\
1 & \ \ \mbox{if } n=-1,0.
\end{array} \right. \]
Finally, once we show the following lemma we will be in position to
prove \teoc{teo2}.

\begin{lem}\label{lem3}
Following the notation in the statement of \teoc{teo2}, let us
write $X_{\varepsilon}-X_0=\lambda X_0+\mu U_0$ and, setting
$\lambda=\sum_{i=1}^{\infty}\lambda_i\varepsilon^i,$ define
 \[
  R(s)=\int_0^{T_0}\left.U_0(\lambda_{k})\right|_{(x,y)=\varphi(t;s)}dt,
 \]
where $\varphi(t;s)$ is the solution of $X_0$ with
$\varphi(0;s)=\xi(s).$ Then $R(s)=\alpha I_0(s)+\beta I_1(s)+\gamma
I_2(s)$ where $(I_0,I_1,I_2)$ is an ECT-system on $\mathcal I$ and
$(\alpha,\beta,\gamma)=\phi(a_k,b_k,c_k,d_k)$ with $\phi$ being a
surjective linear mapping such that
\begin{enumerate}[$(a)$]
\item $\mbox{Ker}(\phi)=\{d_k+3c_k=a_k+3c_k=b_k-c_k=0\}$ in case that $X_0$ is $(S_1^{\ast}),$
\item $\mbox{Ker}(\phi)=\{b_k=c_k=a_k-d_k=0\}$ in case that $X_0$ is $(S_2^{\ast}),$
\item $\mbox{Ker}(\phi)=\{b_k=2a_k+3c_k=2d_k+9c_k=0\}$ in case that $X_0$ is
$(S_3^{\ast})$ or $(\bar S_3^{\ast}).$
\end{enumerate}
\end{lem}

\begin{prova}
First we shall write $R(s)$ as an Abelian integral taking advantage
of the fact that each isochronous center $X_0$ has a first integral
$H,$ see \lemc{lem1}. Thus, if $\gamma_s$ denotes the periodic orbit
of $X_0$ inside the energy level $\{H(x,y)=H\bigl(\xi(s)\bigr)\}$
and $X_0=P_0\partial_x+Q_0\partial_y,$ then
 \[
  R(s)=\int_{\gamma_s}\frac{<\!\nabla\lambda_{k}(x,y),U_0(x,y)\!>}{P_0(x,y)}\,dx,
 \]
where $<\,,\,>$ denotes the inner product. In fact, setting
$\eta(s)\!:=H\bigl(\xi(s)\bigr),$ for simplicity in the exposition
we shall study
 \[
  I(h)\!:=R\bigl(\eta^{-1}(h)\bigr)=\int_{\gamma_h}\frac{<\!\nabla\lambda_{k}(x,y),U_0(x,y)\!>}{P_0(x,y)}\,dx,
 \]
where $\gamma_h$ is now the oval inside the level curve
$\{H(x,y)=h\}.$ Note moreover that, by definition,
$Y_{\varepsilon}=X_{\varepsilon}-X_0$ and
$X_{\varepsilon}-X_0=\lambda X_0+\mu U_0,$ so that
 \begin{equation*}
  \lambda=\frac{Y_{\varepsilon}\wedge U_0}{X_0\wedge U_0}.
 \end{equation*}
We are now in position to prove the result for each different case.

\medskip

\noindent\textbf{The case $\mathbf{(S_2^{\ast})}$.} $ $From the
above equality, taking the expression of $X_0$ and $U_0$ given by
\lemc{lem1} into account, an easy computation shows that
 \[
  \lambda(x,y;\varepsilon)=-{\frac {a(\varepsilon){y}^{2}{x}^{2}+b(\varepsilon){y}^{4}
               -c(\varepsilon){x}^{4}-d(\varepsilon){x}^{2}{y}^{2}}{{x}^{2}+{y}^{2}}},
 \]
and accordingly, since
$\lambda=\sum_{i=1}^{\infty}\lambda_i\varepsilon^i,$ it turns out
that
 \[
  \lambda_{k}(x,y)=-{\frac {a_{k}{y}^{2}{x}^{2}+b_{k}{y}^{4}
               -c_{k}{x}^{4}-d_{k}{x}^{2}{y}^{2}}{{x}^{2}+{y}^{2}}}.
 \]
$ $From now on, when there is no risk of ambiguity, we shall omit
the subscript $k$ for the sake of shortness. In fact it does not
play any role at all because
$\lambda_k(x,y)=F(a_k,b_k,c_k,d_k,x,y)$ with $F$ not depending on
$k.$ Since $H(x,y)=\frac{x^2+y^2}{1-x^2},$ the projection of
$\mathcal P$ on the $x$-axis is $(-1,1)$ and the energy level in
$\gamma_h$ ranges from $h=0$ to $h_0=+\infty.$ It turns out that
 \[
  I(h)=2\int_{\gamma_h}\frac{by^4+(a-d)x^2y^2-cx^4}{y(x^2+y^2)}\, dx=
  \frac{2}{h}\int_{\gamma_h}
  \frac{by^4+(a-d)x^2y^2-cx^4}{y(1-x^2)}\, dx,
 \]
where in the last equality we use that $x^2+y^2=(1-x^2)h.$
Consequently
 \[
  I(h)=\frac{2}{h}\bigl(bI_0(h)+(a-d)I_1(h)-cI_2(h)\bigr)
 \]
with
 \[
  I_0(h)=\int_{\gamma_h}\frac{y^3dx}{1-x^2},\;
  I_1(h)=\int_{\gamma_h}\frac{x^2ydx}{1-x^2}\,\mbox{ and }
  I_2(h)=\int_{\gamma_h}\frac{x^4dx}{y(1-x^2)}.
 \]
Clearly, since $\eta(s)=H\bigl(\xi(s)\bigr)$ is a diffeomorphism and
$\eta(\mathcal I)=(0,+\infty),$ setting $\phi(a,b,c,d)=(b,a-d,-c),$
the proof of \lemc{lem3} will follow for this case once we show that
$(I_0,I_1,I_2)$ is an ECT-system on $(0,+\infty).$ To this end we
shall apply \propc{prop2}. With this aim in view, the application of
\lemc{puja} to
$I_1(h)$ with $k=3$ 
gives
 \[
  I_1(h)=\int_{\gamma_h}f_1(x)y^3dx\,\mbox{ where }f_1(x)=\frac{1-4x^2}{3(1-x^2)}.
 \]
Similarly, but now applying twice \lemc{puja}, first with $k=1$ and
then with $k=3,$ we obtain
 \[
  I_2(h)=\int_{\gamma_h}f_2(x)y^3dx\mbox{ where }f_2(x)=\frac{8x^4-8x^2+1}{1-x^2}.
 \]
On the other hand, note that
 \[
  I_0(h)=\int_{\gamma_h}f_0(x)y^3dx\mbox{ where }f_0(x)=\frac{1}{1-x^2}.
 \]
Our goal with these manipulations is twofold. Firstly, we want
that $y$ has the same exponent in all the Abelian integrals.
Secondly, that this exponent is large enough so that, with the
notation in \propc{prop2}, $m\geqslant n-1.$ Now, by applying this
result, if $(f_0,f_1,f_2)$ is a CT-system on $(0,1),$ then
$(I_0,I_1,I_2)$ is an ECT-system on $(0,+\infty),$ and we are
done. This is clear because $f_0,$ $f_1$ and $f_2$ share the same
denominator and each numerator is an even polynomial of degree
exactly $2i$ for i $=0,1,2.$

\medskip

\noindent\textbf{The case $\mathbf{(S_3^{\ast})}$.} We omit many of
the explanations for the sake of brevity because the proof in this
case is completely analogous to the previous one. Now the Abelian
integral $I(h)=R\bigl(\eta^{-1}(h)\bigr)$ is given by
\[
 I(h)=\frac{2}{h}\int_{\gamma_h}\left(
  \frac{cx^4(1-2x^2)}{(1-3x^2)^3y}+\frac{x^2(1-2x^2)(d-a+4x^2(3a-d))y}{(1-3x^2)^3}
  -\frac{b(1-18x^2+48x^4)y^3}{(1-3x^2)^3}\right)dx,
\]
where we used the transversal commutator and the first integral
given by \lemc{lem1}. The projection of $\mathcal P$ on the $x$-axis
is $\bigl(-\frac{1}{\sqrt{3}},\frac{1}{\sqrt{3}}\bigr)$ and
$h_0=+\infty.$ By applying \lemc{puja} exactly as in the previous
case we obtain
 \[
  I(h)=\frac{2}{3h}\Bigl((d-a-3b+3c)I_0(h)+6(9a+9b-9c-5d)I_1(h)-48(6a+3b-2d)I_2(h)\Bigr),
 \]
where
 \[
  I_i(h)=\int_{\gamma_h}f_i(x)y^3dx\,\mbox{ with
  }f_i(x)=\frac{x^{2i}}{(1-3x^2)^3}\,
  \mbox{ for $i=0,1,2$}.
 \]
By applying \propc{prop2}, $(I_0,I_1,I_2)$ is an ECT-system on
$(0,+\infty)$ because it is clear that $(f_0,f_1,f_2)$ is a
CT-system on $\bigl(0,\frac{1}{\sqrt{3}}\bigr).$ Thus, taking
$\phi(a,b,c,d)=(d-a-3b+3c,9a+9b-9c-5d,6a+3b-2d),$ the result follows
in this case.

\medskip

\noindent\textbf{The case $\mathbf{(\bar S_3^{\ast})}$.} Exactly in
the same way as in the previous cases, the Abelian integral
$I(h)=R\bigl(\eta^{-1}(h)\bigr)$ is given by
\[
 I(h)=\frac{2}{h}\int_{\gamma_h}\left(
  \frac{cx^4(1+2x^2)}{(1+3x^2)^3y}+\frac{x^2(1+2x^2)(d-a+4x^2(d-3a))y}{(1+3x^2)^3}
  -\frac{b(1+18x^2+48x^4)y^3}{(1+3x^2)^3}\right)dx.
\]
The projection of the period annulus $\mathcal{P}$ on the $x$-axis
is $(-\infty,+\infty)$ and $h_0=\frac{4}{27}.$ By applying
\lemc{puja},
 \[
  I(h)=\frac{2}{3h}\Bigl((d-a-3b+3c)I_0(h)-6(9a+9b-9c-5d)I_1(h)-48(6a+3b-2d)I_2(h)\Bigr),
 \]
where
 \[
  I_i(h)=\int_{\gamma_h}f_i(x)y^3dx\,\mbox{ with
  }f_i(x)=\frac{x^{2i}}{(1+3x^2)^3}\,
  \mbox{ for $i=0,1,2$}.
 \]
We are under the hypothesis of Proposition \ref{prop2} and, since it
is obvious that $(f_0,f_1,f_2)$ is a CT-system on $(0,+\infty)$, we
conclude that $(I_0,I_1,I_2)$ is an ECT-system on
$(0,\frac{4}{27})$. Consequently the result follows taking, as
before, $\phi(a,b,c,d)=(d-a-3b+3c,9a+9b-9c-5d,6a+3b-2d).$

\medskip

\noindent\textbf{The case $\mathbf{(S_1^{\ast})}$.} The proof for
the perturbation of this isochronous center is longer and more
complicated than the others because the condition for the first
integral to take the form $A(x)+B(x)y^{2}$, as established in
\propc{prop2}, is not verified. Then we must apply \lemc{lem2}, that
characterizes the Chebyshev property in terms of Wronskians. For
this reason, instead of introducing Abelian integrals, we keep the
expression of $R(s)$ in terms of the solution of the unperturbed
system and, for the sake of convenience, we use complex notation.
Recall that we now study the unfolding
 \[
 \left\{
   \begin{array}{l}
    \dot x=-y+\bigl(-3+a(\varepsilon)\bigr)x^2y+\big(1+b(\varepsilon)\bigr)y^3, \\[3pt]
    \dot{y}=x + \bigl(1+c(\varepsilon)\bigr)x^3+\bigl(-3+d(\varepsilon)\bigr)xy^2.
   \end{array}
  \right.
 \]
Taking the commutator given by \lemc{lem1}, an easy computation
shows that
 \[
  \lambda(x,y;\varepsilon)=\frac{x^2(1+x^2-3y^2)\bigl(c(\varepsilon)x^2+d(\varepsilon)y^2\bigr)
    -y^2(1+3x^2-y^2)\bigl(a(\varepsilon)x^2+b(\varepsilon)y^2\bigr)}
    {(x^2+y^2)\bigl(1+2(x^2-y^2)+(x^2+y^2)^2\bigr)}.
 \]
The coordinate transformation $z=x+iy$ brings the unfolding to $\dot
z=f(z)+ip(z,\bar z)$ with
 \[
  f(z)=iz(1+z^2)\quad \mbox{ and } \quad p(z,\bar z)=(\alpha  z^3+\beta z^2\bar{z}+\gamma z\bar{z}^2+
  \delta \bar{z}^3)/8,
 \]
where $\alpha=b+c-d-a,$ $\beta=3c+d-a-3b,$ $\gamma=a+3b+3c+d$ and
$\delta=a-b+c-d.$ An easy computation shows that the solution of the
unperturbed system, i.e., $\dot z=f(z),$ with initial condition
$z=h\in (0,+\infty)$ at $t=0$ is given by
\[
 \varphi(t;h)=\frac{he^{it}}{\sqrt{1+h^2-e^{2it}h^2}}.
\]
Accordingly, following the usual notation
$\alpha(\varepsilon)=\sum_{i=1}^{\infty}\alpha_i\varepsilon^i,$ we
split the function under consideration as
 \[
  R(h)=\int_0^{2\pi}\!\left.<\!\nabla\lambda_{k}(x,y),U_0(x,y)\!>\right|_{x+iy=\varphi(t;h)}dt=
  \alpha_{k} I_3(h)+\beta_{k} I_2(h)+\gamma_{k} I_1(h)+\delta_{k} I_0(h).
 \]
Long but straightforward manipulations show that $I_3\equiv 0$ and
that
 \[
 R(h)=\frac{-4h^4}{(1+h^2)^2}\Bigl(\beta_{k}\bar I_2(h)+\gamma_{k}\pi+\delta_{k}\bar
 I_0(h)\Bigr)
 \]
where, setting $\mu(h)=\frac{2h\sqrt{1+h^2}}{1+2h^2},$
\begin{align*}
 &\bar I_2(h)\!:=\frac{2h^2}{1+2h^2} \,
 \mathcal{K}\bigl(\mu(h)\bigr)-\frac{2+4h^2}{h^2} \, \mathcal{E}\bigl(\mu(h)\bigr)
 \intertext{and}
 &\bar I_0(h)\!:=\frac{(1+2h^2)(1+2h^2+2h^4)}{h^6}\,\mathcal{E}\bigl(\mu(h)\bigr)
   -\frac{(1+h^2+h^4)(1+3h^2+3h^4)}{h^6(1+2h^2)}\,\mathcal{K}\bigl(\mu(h)\bigr).
\end{align*}

The fact that $I_3\equiv 0$ is not unexpected at all. If
$\beta=\gamma=\delta=0,$ then the perturbed system is holomorphic,
it does not depend on $\bar z,$ and consequently the center is
isochronous for all $\varepsilon,$ so that
$T'(s;\varepsilon)\equiv 0$, see \cite{Villarini}.

At this point, taking $\phi(a,b,c,d)=(3c+d-a-3b,a+3b+3c+d,a-b+c-d),$
the proof of \lemc{lem3} for the case $(S_1^{\ast})$ reduces to the
verification of the fact that $(\pi,\bar I_2,\bar I_0)$ is an
ECT-system on $(0,+\infty).$ Since the first function is a non-zero
constant, we compute the two-dimensional Wronskian by using
\lemc{lem5},
\[
 W[\pi,\bar I_2](h)=\bar I_2'(h)=\frac{2+2h^2}{h^3}
 \left( \mathcal{E}\bigl(\mu(h)\bigr)+\frac{\mathcal{K}\bigl(\mu(h)\bigr)}{1+2h^2}\right).
\]
Clearly it is different from zero for all $h\in (0,+\infty)$ because $\mathcal{E}$ and
$\mathcal{K}$ are strictly positive functions. Next, taking \lemc{lem5} into account
again, we compute the three-dimensional Wronskian. The key point is that it factorizes as
\[
 W\bigl[\pi,\bar I_2,\bar I_0\bigr](h)=\left|
\begin{array}{cc}
\bar I_2' & \bar I_0' \\[4pt] \bar I_2'' & \bar I_0''
\end{array}\right|=
\frac{72(1+h^2)^3}{h^{11}}\Bigl(\mathcal{E}\bigl(\mu(h)\bigr)+L_{+}(h)\,\mathcal{K}\bigl(\mu(h)\bigr)\Bigr)
\Bigl(\mathcal{E}\bigl(\mu(h)\bigr)+L_{-}(h)\,\mathcal{K}\bigl(\mu(h)\bigr)\Bigr),
\]
where
\[
 L_{\pm}(h)\!:=\frac{-1-2h^2-2h^4\pm 2(1+4h^2+5h^4+2h^6+h^8)^{1/2}}{3(1+2h^2)^2} \, .
\]
Note that $L_-(h)<0$ for all $h\in (0,+\infty).$ It is also easy to show that $L_+(h)>0$
for all $h\in (0,+\infty).$ Therefore, we will see that this Wronskian does not vanish
once we prove that
 \[
  \mathcal R(h)\!:=\mathcal{E}\bigl(\mu(h)\bigr)+L_{-}(h)\,\mathcal{K}\bigl(\mu(h)\bigr)
    \neq 0\,\mbox{ for all $h\in (0,+\infty).$}
 \]
$ $From now on we  will use the variable $u=\mu(h)$ because then
the expressions that we obtain are shorter. Thus, due to
$\mu(0,+\infty)=(0,1),$ we must show that
 \[
  \mathcal L(u)\!:=\mathcal R\bigl(\mu^{-1}(u)\bigr)=\mathcal{E}(u)+\frac{1}{6}\left(
   u^2-2 -\sqrt{16-16u^2+u^4}\,\right) \mathcal{K}(u)
 \]
does not vanish for all $u\in (0,1).$  To prove this claim we first
note that by applying \lemc{lem5} one can check that $\mathcal L(u)$
verifies the differential equation
 \begin{equation}\label{ed}
  \mathcal L''(u)=g_1(u)\mathcal L'(u)+g_0(u)\mathcal L(u),
 \end{equation}
where
\begin{align*}
&g_1(u)=\frac{48-64u^2+17u^4}{u(1-u^2)(16-16u^2+u^4)}+
\frac{\sqrt{16-16u^2+u^4}}{u(1-u^2)}\intertext{and}
&g_0(u)=\frac{16-12u^2+
(8-u^2)\sqrt{16-16u^2+u^4}}{(1-u^2)(16-16u^2+u^4)}.
\end{align*}
On the other hand, from \lemc{lem5} once again, we get that
$\mathcal L(u)=-\frac{3\pi}{4096}u^8+\op(u^{9}).$ We are now in
position to prove the claim. By contradiction, assume that there
exists some $u_1\in (0,1)$ such that $\mathcal L(u_1)=0.$ Then,
since $\mathcal L(u)<0$ for $u\approx 0,$ this forces the existence
of a local minimum of~$\mathcal L,$ say $u_0,$ with $\mathcal
L(u_0)<0.$ The evaluation of the differential equation in \refc{ed}
at $u=u_0$ shows that $\mathcal L''(u_0)=g_0(u_0)\mathcal L(u_0).$
Since $g_0(u)>0$ for all $u\in (0,1),$ this implies that $\mathcal
L''(u_0)<0,$ which contradicts the fact that $u=u_0$ is a local
minimum. Hence the claim is true and, therefore, the
three-dimensional Wronskian is different from zero. In short,
$(\pi,\bar I_2,\bar I_0)$ is an ECT-system on $(0,+\infty)$ and this
completes the proof of the result for the case $(S_1^{\ast}).$
\end{prova}

\begin{prooftext}{Proof of \teoc{teo2}}
We prove the result for the perturbation of the isochronous center
$(S_1^{\ast})$ only because the other cases follow exactly in the
same way. Thus, consider the vector field $X_0$ with the isochronous
center $(S_1^{\ast})$ at the origin and, see \lemc{lem1}, let $U_0$
be its commutator. In addition, let $\Phi$ be a linearization of
$X_0.$ As usual we consider the family of centers
$X_{\varepsilon}=X_0+Y_{\varepsilon}$ and we denote by
$T(s;\varepsilon)$ the period function of $X_{\varepsilon}$ using a
solution of $U_0$ as transversal section.

We claim that if $T_0'\equiv T_1'\equiv\ldots\equiv T_k'\equiv 0,$
then there exists an analytic family of diffeomorphisms
$\{\Phi_{\varepsilon}^k\},$ in a neighbourhood of the origin,  such
that~$\Phi_{\varepsilon}^k$ linearizes
$j^k\bigl(X_{\varepsilon}\bigr)$ for all $\varepsilon\approx 0.$
The proof follows by induction on $k.$ The case $k=0$ is trivial
because $j^{\,0}(X_{\varepsilon})=X_0$ and so we can take
$\Phi_{\varepsilon}^0=\Phi.$ Assume that the claim is true for
$k=n$ and let us show its validity for $k=n+1.$ So suppose that
$T_0'\equiv T_1'\equiv\ldots\equiv T_n'\equiv T_{n+1}'\equiv 0$
and, by the induction hypothesis, that there exists a
linearization $\Phi_{\varepsilon}^k$ of $j^{k}(X_{\varepsilon})$
for $k=0,1,\ldots,n.$ Then, by applying $n+1$ times \teoc{teo1},
 \[
  T_k'(s)=-\int_0^{T_0}\left.U_0(\lambda_k)\right|_{(x,y)=\varphi(t;s)}dt\,
  \mbox{ for $k=1,2,\ldots,n+1,$}
 \]
and consequently, from $(a)$ in \lemc{lem3},
 \[
  d_k+3c_k=a_k+3c_k=b_k-c_k=0\,\mbox{ for $k=1,2,\ldots,n+1.$}
 \]
(This follows from the fact that each $T_k'$ is a linear combination
of three functions forming an ECT-system with the coefficients
vanishing simultaneously only in case that the above relations
hold.) Therefore it turns out that we can write
 \[
  j^{\,n+1}(X_{\varepsilon})=\bigl(-y-3\kappa(\varepsilon)x^2y +\kappa(\varepsilon)y^3\bigr)
   \partial_x+\bigl(x+\kappa(\varepsilon)x^3-3\kappa(\varepsilon)xy^2\bigr)\partial_y
 \]
with $\kappa(\varepsilon)\!:=1+\sum_{k=1}^{n+1}c_i\varepsilon^k.$
Now, if we define
$\Psi_{\varepsilon}(x,y)=\bigl(\kappa(\varepsilon)\,x,\kappa(\varepsilon)\,y\bigr),$
then one can easily verify that it holds
$\left(\Psi_{\varepsilon}\right)^{\ast}\left(j^{\,n+1}(X_{\varepsilon})\right)=X_0.$
Accordingly
$\Phi_{\varepsilon}^{n+1}\!:=\Phi\circ\Psi_{\varepsilon}$ provides a
linearization of the $(n+1)$-jet of $X_{\varepsilon}.$ Thus the
claim is true.

It is clear that the result under consideration only makes sense in
case that the perturbation is not isochronous, so there exists some
$\ell\geqslant 0$ such that $T_0'\equiv T_1'\equiv\ldots\equiv
T_\ell'\equiv 0$ and $T_{\ell+1}'\not\equiv 0.$ Then, on account of
the claim, there exists an analytic family of linearizations
$\{\Phi_{\varepsilon}^\ell\}$ of $j^{\ell}(X_{\varepsilon})$ and so,
by applying \teoc{teo1},
 \[
  T_{\ell+1}'(s)=-\int_0^{T_0}\left.U_0(\lambda_{\ell+1})\right|_{(x,y)=\varphi(t;s)}dt\,
  \mbox{ for all $s\in\mathcal I,$}
 \]
Hence, from \lemc{lem3}, $T_{\ell+1}'(s)=\alpha I_0(s)+\beta
I_1(s)+\gamma I_2(s)$ where $(I_0,I_1,I_2)$ is an ECT-system on
$\mathcal I$ and
$(\alpha,\beta,\gamma)=\phi(a_{\ell+1},b_{\ell+1},c_{\ell+1},d_{\ell+1})$
for some surjective linear mapping. Consequently, since
$T_{\ell+1}'$ is not identically zero, it can have at most two zeros
on $\mathcal I.$ Moreover, since $\phi$ is exhaustive, one can
choose $a_{\ell+1},b_{\ell+1},c_{\ell+1}$ and $d_{\ell+1}$ such that
$T_{\ell+1}'(s)=0$ has exactly $0,$ $1$ or $2$ roots for
$s\in\mathcal I$ (recall \obsc{hison}). This completes the proof of
the result.
\end{prooftext}

\section{Proof of \teoc{Loud}}\label{sect4}

Recall that we consider the unfolding $X_{\eps}=X_{0}+Y_{\eps},$
where $X_0$ is the vector field associated to each one of the
isochrones in \refc{llista2} and $Y_{\eps}$ is the perturbation in
\refc{pert2}. Note first that, by means of the transformation
$(x,y,t)\longmapsto (\eta x,\eta y,\eta t)$ with $\eta=1+b(\eps),$
there is no loss of generality in assuming $b\equiv 0,$ i.e., that
the unfolding is given by
 \begin{equation}\label{pert3}
  \sist{-y+xy,}{x+\bigl(d_0+d(\eps)\bigr)x^2+\bigl(f_0+f(\eps)\bigr)y^2,}
 \end{equation}
where $d$ and $f$ are analytic functions vanishing at $\eps=0,$
and the value of $(d_0,f_0)$ is $(-\frac{1}{2},\frac{1}{2}),$
$(0,1),$ $(1,\frac{1}{4})$ and $(-\frac{1}{2},2)$ in case that we
perturb the isochronous center $(S_1),$ $(S_2),$ $(S_3)$ and
$(S_4)$, respectively.

In order to show \teoc{Loud} we must first take, for each
isochronous center, a commutator $U_0$ of $X_0$ defined in the
whole period annulus. For instance we can use the ones in
\cite{ChaSab,MarMosRou} but, as it will be clear in a moment, we
do not need their concrete expression. Then, as we did in the
previous section for the Pleshkan's isochrones, we decompose the
perturbation as $X_{\varepsilon}-X_0=\lambda X_0+\mu U_0.$ The
next result is the counterpart of \lemc{lem3} and it follows after
translating the results in \cite{GasYu} to the language that we
use here.
\begin{lem}[Gasull-Yu]\label{LemGas}
Setting $\lambda=\sum_{i=1}^{\infty}\lambda_i\varepsilon^i,$ define
 \[
  R(s)=\int_0^{T_0}\left.U_0(\lambda_{k})\right|_{(x,y)=\varphi(t;s)}dt,
 \]
where $\varphi(t;s)$ is the solution of $X_0$ such that
$\varphi(0;s)=\xi(s)$ and \map{\xi}{\mathcal I}{\R^2} is a fixed
solution of $U_0.$ Then
 \[
  R(s)=d_k I_0(s)+f_k I_1(s)\,\mbox{ for all $s\in\mathcal I,$}
 \]
where:
 \begin{enumerate}[$(a)$]
  \item $(I_0,I_1)$ is an ECT-system on $\mathcal I$ in case that
        $(d_0,f_0)\in\bigl\{(0,1),(0,\frac{1}{4}),(-\frac{1}{2},2)\bigr\}.$
  \item $I_0(s)=I_1(s)\neq 0$ for all $s\in\mathcal I$ in case that $(d_0,f_0)=(-\frac{1}{2},\frac{1}{2}).$
 \end{enumerate}
\end{lem}

This result follows from the proof of Theorem 3 in \cite{GasYu}, but
it is worth making some comments. Recall that to parameterize the
periodic orbits we take a commutator $U_0$ of $X_0,$ i.e., such that
$[X_0,U_0]=0,$ analytic in $\mathcal P \cup \{(0,0)\}$. This
provides a parameterization of the period function of $X_{\eps},$
say $T(s;\eps)=2\pi+T_1(s)\eps+\op(\eps),$ for which
 \[
  T_1'(s)=\int_0^{T_0}\left.U_0(\lambda_1)\right|_{(x,y)=\varphi(t;s)}dt.
 \]
(This expression is well known but it can be viewed as the case
$k=0$ of \teoc{teo1}.) Instead of a commutator, the authors in
\cite{GasYu} use a vector field $\widehat U_0$ such that
$[X_0,\widehat U_0]=\beta\,\widehat U_0$ for some function
$\beta.$ This yields of course to another parameterization of the
period function of $X_{\eps},$ say $\widehat
T(s;\eps)=2\pi+\widehat T_1(s)\eps+\op(\eps).$ The expression of
$\widehat T_1'(s)$ is slightly different from the previous one,
but we still can take advantage of their result. Indeed, it is
clear that there exists a diffeomorphism $\zeta$ verifying
$T(s;\eps)=\widehat T\bigl(\zeta(s);\eps\bigr),$ so that
$T_1'(s)=\zeta'(s)\widehat T_1'\bigl(\zeta(s)\bigr)$, see
\obsc{indep}. The proof of Theorem 3 in \cite{GasYu} shows that
$\widehat T_1'(s)=d_1\hat I_0(s)+f_1\hat I_1(s),$ where $\hat I_0$
and $\hat I_1$ are functions verifying, for each case, statements
$(a)$ and $(b)$ in \lemc{LemGas}. This is the key point because
then
 \[
 \int_0^{T_0}\left.U_0(\lambda_1)\right|_{(x,y)=\varphi(t;s)}dt
 =d_1 I_0(s)+f_1 I_1(s),\mbox{ with $I_i(s)\!:=\zeta'(s)\hat I_i\bigl(\zeta(s)\bigr),$}
 \]
and the result follows due to $\zeta'(s)\neq 0$ for all $s.$ (To
be more precise, this proves the case $k=1,$ however the subindex
does not play any role at all because
$\lambda_k(x,y)=F(d_k,f_k,x,y)$ with $F$ not depending on $k.$)

Let us note that analogous computations to the ones carried out in
the previous section and the application of the criterion given in
\cite{GMV} also yield to the proof of \lemc{LemGas}. For the sake of
shortness we prefer to take advantage of the results in \cite{GasYu}
instead.

\begin{prooftext}{Proof of \teoc{Loud}}
Let us fix any $(d_0,f_0)\in
\bigl\{(0,1),(0,\frac{1}{4}),(-\frac{1}{2},2)\bigr\}$ and consider
the unfolding $X_{\eps}$ in~\refc{pert3} of the isochronous center
$X_0.$ Let us also take the corresponding commutator $U_0$ of $X_0$
and, as usual, denote the period function of $X_{\varepsilon}$ using
a solution of $U_0$ as transversal section by $T(s;\varepsilon)$.
Finally, let $\Phi$ be a linearization of $X_0.$

We claim that if $T_0'\equiv T_1'\equiv\ldots\equiv T_k'\equiv 0,$
then $j^k\bigl(X_{\varepsilon}\bigr)=X_0.$ The proof follows by
induction on $k.$ The case $k=0$ is trivial because
$j^{\,0}(X_{\varepsilon})=X_0.$ Assume that the claim is true for
$k=n$ and let us show its validity for $k=n+1.$ So suppose that
$T_0'\equiv T_1'\equiv\ldots\equiv T_{n+1}'\equiv 0$ and, by the
induction hypothesis, that $j^n(X_{\varepsilon})=X_0.$ Then $\Phi$
is a linearization of $j^k(X_{\varepsilon})=X_0$ for
$k=0,1,\ldots,n$ and, by applying $n+1$ times \teoc{teo1}, we have
that
 \[
  T_k'(s)=-\int_0^{T_0}\left.U_0(\lambda_k)\right|_{(x,y)=\varphi(t;s)}dt\,
  \mbox{ for $k=1,2,\ldots,n+1.$}
 \]
Hence, thanks to $(a)$ in \lemc{LemGas}, $T_1'\equiv
T_2'\equiv\ldots\equiv T_{n+1}'\equiv 0$ implies $d_k=f_k=0$ for
$k=1,2,\ldots,n+1,$ so that  $j^{\,n+1}(X_{\varepsilon})=X_0.$ This
proves the claim.

It is clear that the result under consideration only makes sense in
case that the perturbation is not isochronous, so there exists some
$\ell\geqslant 0$ such that $T_0'\equiv T_1'\equiv\ldots\equiv
T_\ell'\equiv 0$ and $T_{\ell+1}'\not\equiv 0.$ Then, on account of
the claim, $\Phi$ is a linearization of
$j^{\ell}(X_{\varepsilon})=X_0$ and, by applying \teoc{teo1},
 \[
  T_{\ell+1}'(s)=-\int_0^{T_0}\left.U_0(\lambda_{\ell+1})\right|_{(x,y)=\varphi(t;s)}dt\,
  \mbox{ for all $s\in\mathcal I.$}
 \]
Thus, from \lemc{lem3}, $T_{\ell+1}'(s)=d_{\ell+1} I_0(s)+f_{\ell+1}
I_1(s)$ where $(I_0,I_1)$ is an ECT-system on $\mathcal I.$
Consequently, since $T_{\ell+1}'$ is not identically zero, it can
have at most one zero for $s\in\mathcal I.$ Moreover, on account of
\obsc{hison}, one can choose $d_{\ell+1}$ and $f_{\ell+1}$ such that
$T_{\ell+1}'(s)=0$ has exactly $0$ or $1$ root for $s\in\mathcal I.$
This completes the proof of the result.
\end{prooftext}

\begin{figure}[t]
  \begin{center}
   \epsfig{file=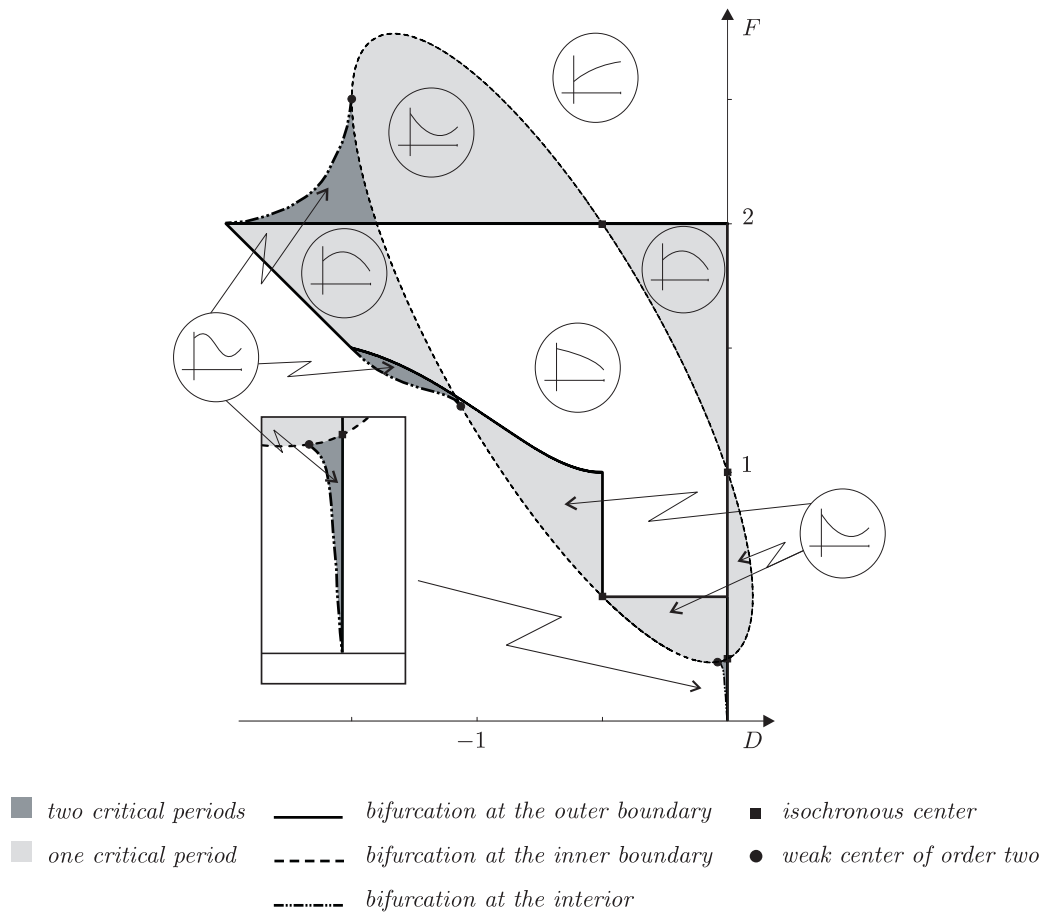,scale=0.9}
  \end{center}
  \caption{Conjectural bifurcation diagram of the period function of the dehomogenized Loud's systems.}
  \label{dib}
 \end{figure}

We conclude the paper with some final remarks about how \teoc{Loud}
and \figc{dib} fit together.
This figure gathers some of the results and conjectures that appear
in \cite{MMV2} on the bifurcation diagram of the period function of
the so-called dehomogenized Loud's systems
 \[
  \sist{-y+xy,}{x+Dx^2+Fy^2.}
 \]
It may seem surprising that, thanks to \teoc{teo1}, it is possible
to obtain a result for arbitrary perturbations by taking advantage
of a result, namely the one in \cite{GasYu}, that holds only for
perturbations that are linear in~$\eps$. Roughly speaking, this
shows that the linear perturbations already give all the critical
periods that you can get with an arbitrary perturbation. Certainly
this is not a general property but it explains why we could not
obtain the result for the perturbation of the isochronous center
$(d_0,f_0)=(-\frac{1}{2},\frac{1}{2}),$ not even for the linear
ones. Indeed, a linear perturbation corresponds in the parameter
plane to move slightly travelling on a straight line that passes
through the isochronous center. In addition if you can get
$k\geqslant 0$ critical periods moving in one direction, then you
also get $k$ in the opposite one (changing the signum of $\eps$).
Clearly, see \figc{dib}, this is what it happens with the
isochronous centers $(d_0,f_0)\in \bigl\{(0,1),
(-\frac{1}{2},2),(0,\frac{1}{4})\bigr\}$ but not with
$(d_0,f_0)=(-\frac{1}{2},\frac{1}{2}).$ Let us note moreover that
if the conjectural diagram in \figc{dib} is true, then the
bifurcation of $(d_0,f_0)=(0,\frac{1}{4})$ giving rise to a
critical period arising from $\mathcal P$ also causes the
emergence of another one from $\partial\mathcal P$ at the same
time. The bifurcation of the latter critical period would be the
counterpart, for the period function, of the bifurcation of an
\emph{alien limit cycle} as defined in \cite{CDR05}. (In that
paper, the authors present a way to study the perturbations from a
Hamiltonian 2-saddle cycle which can produce limit cycles that
cannot be detected using zeroes of the Abelian integral, even when
it is generic.) On account of \teoc{Loud}, this, let us say,
\emph{alien critical period} does not come from a zero of the
Abelian integral related with the derivative of the period
function.

On the other hand, Chicone and Jacobs prove in \cite{ChiJac} that
there are perturbations of each one of the four isochrones giving
rise to exactly one critical period bifurcating from the inner
boundary of $\mathcal P$. \teoc{Loud} shows that there are
perturbations of the isochronous centers $(d_0,f_0)\in
\bigl\{(0,1), (-\frac{1}{2},2),(0,\frac{1}{4})\bigr\}$ giving rise
to exactly one critical period bifurcating from $\mathcal P.$
These two facts do not contradict the conjecture in \figc{dib}
near $(d_0,f_0)=(0,1)$ and $(d_0,f_0)=(-\frac{1}{2},2),$ it simply
shows that both perturbations are different.

\bibliographystyle{plain}

\begin{thebibliography}{99}

\bibitem{Abram} M. Abramowitz and I.A. Stegun, ``Handbook
of mathematical functions with formulas, graphs, and mathematical
tables'' Reprint of the 1972 edition. Dover Publications, Inc., New
York, 1992.

\bibitem{CDR05} M. Caubergh, F. Dumortier and R. Roussarie, {\it Alien limit cycles near a Hamiltonian 2-saddle
cycle}, C. R. Acad. Sci. Paris, Ser. I {\bf 340} (2005) 587–-592.

\bibitem{ChiJac} C. Chicone and M. Jacobs, {\it Bifurcation of critical periods
for plane vector fields,} Trans. Amer. Math. Soc. \textbf{312}
(1989) 433--486.

\bibitem{CimGasSil} A. Cima, A. Gasull and P.R. da Silva, {\it On the number of critical periods for planar
polynomial systems,} Nonlinear Anal. {\bf 69} (2008) 1889–-1903.

\bibitem{ChaSab} J. Chavarriga and M. Sabatini, {\it A survey of
isochronous centers,}  Qual. Theory Dyn. Syst. {\bf 1} (1999) 1--70.

\bibitem{Dulac} H. Dulac, {\it D\'etermination et
int\'egration d'une certaine classe d'\'equations
diff\'erentielles ayant pour point singulier un centre,} Bull.
Soc. Math. France S\'er. (2) {\bf 32} (1908), 230--252.

\bibitem{Fran} J.-P. Fran\c{c}oise, {\it The successive derivatives of the period function
of a plane vector field}, J. Differential Equations {\bf 146}
(1998) 320--335.

\bibitem{FreGasGui1} E. Freire, A. Gasull and A. Guillamon, {\it Period function for
perturbed isochronous centres}, Qual. Theory Dyn. Syst. {\bf 3}
(2002) 275–-284.

\bibitem{FreGasGui2} E. Freire, A. Gasull and A. Guillamon, {\it First derivative of
the period function with applications}, J. Differential Equations
{\bf 204} (2004) 139--162.

\bibitem{GasYu} A. Gasull and Jiang Yu, {\it On the critical periods of perturbed isochronous
centers,} J. Differential Equations {\bf 244} (2008) 696-–715.

\bibitem{GasZhao} A. Gasull and Yulin Zhao, {\it Bifurcation of critical periods
from the rigid quadratic isochronous vector field,}  Bull. Sci.
Math. {\bf 132} (2008) 292–-312

\bibitem{GMV} M. Grau, F. Ma\~{n}osas and J.
Villadelprat, {\it A Chebyshev criterion for Abelian integrals,}
preprint (2008) {\tt arXiv:0805.1140v2 [math.DS]}.

\bibitem{Karlin} S. Karlin and W. Studden, ``Tchebycheff systems:
with applications in analysis and statistics'', Interscience
Publishers, 1966.

\bibitem{KMS} I. Kol\'{a}\v r, P. Michor and J. Slov\'{a}k, ``Natural operations in differential
geometry'', Springer--Verlag, Berlin, 1993.

\bibitem{Loud} W.S. Loud, {\it Behaviour of the period of solutions of certain plane autonomous systems near
centers,} Contrib. Differential Equations \textbf{3} (1964) 21-–36.

\bibitem{Malkin} K.E. Malkin, {\it Criteria for center for a differential system}, Volzhskii. Matem. Sbornik
{\bf 2} (1964) 87-–91.

\bibitem{Mar} P. Marde\v si\'c, ``Chebyshev systems and the versal unfolding of the cusp of order~$n$'',
Travaux en cours, vol.~57, Hermann, Paris, 1998.

\bibitem{MMV1} P. Mardesic, D. Mar\'{\i}n and J. Villadelprat,
{\it On time function of the Dulac map for families of meromorphic
vector fields,} Nonlinearity \textbf{16} (2003) 855-–881.

\bibitem{MMV2} P. Marde\v si\'c, D. Mar\'{\i}n and J. Villadelprat, {\it The period function of reversible quadratic
centers,}  J.~Differential Equations \textbf{224} (2006) 120--171.

\bibitem{MV} D. Mar\'{\i}n and J. Villadelprat, {\it On the return time function around
monodromic polycycles,} J. Differential Equations \textbf{228}
(2006) 226–-258.

\bibitem{MarMosRou} P. Marde\v{s}i\'c, L. Moser-Jauslin and C.
Rousseau, {\it Darboux linearization and isochronous centers with a
rational first integral}, J. Differential Equations {\bf 134} (1997)
216--268.

\bibitem{Pleshkan} I. Pleshkan, {\it A new method of investigating the isochronicity of a system of two differential
equations}, Differential Equations {\bf 5} (1969) 796-–802.

\bibitem{RouToni} C. Rousseau and B. Toni, {\it Local bifurcation of critical periods in vector fields with
homogeneous nonlinearities of the third degree}, Canad. Math. Bull.
{\bf 36} (1993) 473--484.

\bibitem{Sabatini} M. Sabatini, {\it Characterizing isochronous centres by Lie brackets},
Differential Equations Dynam. Systems {\bf 5} (1997) 91--99.

\bibitem{Villarini} M. Villarini, {\it  Regularity properties of the period function near
a centre of a planar vector field}, Nonlinear Analysis T.M.A. {\bf
19} (1992) 787–-803.

\end{thebibliography}

\end{document}